\documentclass[12pt]{amsproc}
\usepackage[cp1251]{inputenc}
\usepackage[T2A]{fontenc} 
\usepackage[russian]{babel} 
\usepackage{amsbsy,amsmath,amssymb,amsthm,graphicx,color} 
\usepackage{amscd}
\def\le{\leqslant}
\def\ge{\geqslant}

\newtheorem{prop}{Предложение}

\theoremstyle{definition}

\theoremstyle{remark}

\begin {document}
\unitlength=1mm
\title[Числовые полугруппы, суммы Дедекинда]
{Числовые полугруппы, суммы Дедекинда и многочлены Александера торических узлов}
\author{Г. Г. Ильюта}
\email{ilgena@rambler.ru}
\address{}
\thanks{Работа поддержана грантом РФФИ-20-01-00579}


\begin{abstract}
Мы свяжем суммы Дедекинда и многочлены Александера торических узлов.

We connect Dedekind sums and Alexander polynomials of torus knots.
\end{abstract}

\maketitle
\tableofcontents

\section{Введение}

  Мы свяжем суммы Дедекинда с многочленом Александера $(a,b)$-торического узла, который интерпретирован как полугрупповой многочлен числовой полугруппы $S_{a,b}$, порождённой взаимно простыми натуральными числами $a$ и $b$. Возможность такой интерпретации вытекает из формулы для ряда Гильберта числовой полугруппы $S_{a,b}$ \cite{15}, \cite{14}, \cite{12}, \cite{4}
$$
H_{S_{a,b}}(q):=\sum_{k\in S_{a,b}}q^k=\frac{1-q^{ab}}{(1-q^a)(1-q^b)}.    \eqno (1)
$$
Добавим ещё одно доказательство: формула (1) является прямым следствием равенств из \cite{11} и \cite{10}, соответственно,
$$
C_{S_{a,b}}:=\mathbb Z_{\ge 0}\setminus S_{a,b}=\{ab-ia-jb:0<i<b,0<j<a,ia+jb<ab\},
$$
$$
\sum_{\substack{0\le i<b,0\le j<a\\ia+jb<ab}}q^{ia+jb}=\frac{1-q^{ab}}{(1-q^a)(1-q^b)}-\frac{q^{ab}}{1-q} 
$$
и очевидного соотношения
$$
C_{S_{a,b}}(q):=\sum_{k\in C_{S_{a,b}}}q^k=\frac{1}{1-q}-H_{S_{a,b}}(q).        
$$
П. Мори в \cite{12} относит формулу (1) к фольклору, отмечая, что для простых чисел $a$ и $b$ она без доказательства упоминается в \cite{7}. Формула для многочлена Александера $(a,b)$-торического узла
$$
A_{S_{a,b}}(q)=\frac{(1-q^{ab})(1-q)}{(1-q^a)(1-q^b)}
$$
была доказана В. Бурау в 1932 г. Из приведённых равенств вытекают соотношения
$$
A_{S_{a,b}}(q)=(1-q)H_{S_{a,b}}(q)=1-(1-q)C_{S_{a,b}}(q).
$$
Аналогичное равенство имеет место для полугруппового многочлена $A_S(q)$ произвольной числовой полугруппы $S$ \cite{12}.

  Обобщенные суммы Дедекинда обычно определяются либо как тригонометрические суммы, либо как суммы, содержащие функцию целая часть. Будем называть соответствующие равенства формулами тригонометрического или рационального типов. Мы покажем, что в формулах тригонометрического типа для некоторых обобщений сумм Дедекинда естественным образом появляются значения многочлена Александера торического узла в корнях из единицы (или в корнях из единицы, умноженных на переменную). Такая же ситуация возникает в формуле Фокса (и её обобщениях) для порядка группы гомологий циклического накрытия дополнения к узлу в $S^3$, но мы не будем заниматься здесь этой темой. Отметим только, что в нашей терминологии формула Фокса и её обобщения являются формулами тригонометрического типа (циклический результант многочлена Александера и его обобщения) и для алгебраических узлов с помощью соотношений из п.2 эти равенства можно представить в виде формул рационального типа. Аналогично можно преобразовать циклический результант и его обобщения для полугруппового многочлена произвольной числовой полугруппы.

  Идея использования при изучении числовых полугрупп дискретного преобразования Фурье (или, что равносильно, соотношений ортогональности для характеров циклической группы) содержится в \cite{1}. В не совсем явном виде мы используем эту идею при доказательстве Предложения 1.

  В п.2 мы выразим через $C_{S_{a,b}}(q)$ (или, что равносильно, через $A_{S_{a,b}}(q)$) ключевую компоненту $[ak/b]$ ($[x]$ -- целая часть $x$) в определениях нескольких обобщений сумм Дедекинда, а также $q^{b[ak/b]}$, что позволит в п.4 получить содержащую многочлены $C_{S_{a,b}}(\epsilon_b^jq)$ формулу для $c(q^b,t;a,b)$, где
$$
c(q,t;a,b):=\sum_{k=1}^{b-1}q^{[ak/b]}t^{k-1}
$$
 -- многочлены Дедекинда-Карлитца \cite{8}, \cite{6}, \cite{5}, $\epsilon_b=e^{\frac{2\pi i}{b}}$.

  В п.3 мы докажем содержащую числа $C_{S_{a,b}}(\epsilon_b^k)$ и специализации многочленов Мириманова (или многочленов Апостола-Бернулли) формулу тригонометрического типа для сумм 
$$
V_{m,n}(a,b):=\sum_{k=1}^{b-1}k^m\left[\frac{ak}{b}\right]^n.
$$
Эти суммы (мы будем называть их суммами Вороного) появились в 1889 г. в студенческой работе Г. Ф. Вороного при доказательстве сравнения для чисел Бернулли \cite{16}, \cite{2}. А именно, в этой работе представлены в виде линейных комбинаций сумм Вороного степенные суммы $M_k(1,m)$, где $M_k(\lambda,m)$ -- многочлены Мириманова (см. п.3). Мы получим в некотором смысле обращение этой формулы -- представление сумм Вороного в виде линейных комбинаций значений многочленов Мириманова в корнях из единицы. Многочлены Мириманова и многочлены Апостола-Бернулли связаны известной формулой \cite{3}, \cite{17}. Некоторые формулы для сумм Вороного, в частности, теоремы взаимности для малых $m$ и $n$, доказаны в \cite{9}, \cite{6}. Суммы $V_{1,1}(a,b)$ связаны простой формулой с классическими суммами Дедекинда \cite{13} (см. п.5).

  Многочлены Дедекинда-Карлитца и суммы Вороного можно объединить двумя способами, а именно, они находятся среди специализаций двух определённых в п.4 многочленов.

  В п.5 для классических сумм Дедекинда мы докажем формулу тригонометрического типа, содержащую числа $C_{S_{a,b}}(\epsilon_b^k)$.

  Неформальное следствие этой статьи состоит в следующем: два представления сумм Дедекинда -- тригонометрическое и рациональное -- аналогичны двум способам определить числовую полугруппу с помощью конечных множеств -- через дополнение в $\mathbb Z_{\ge 0}$ или с помощью множеств Апери. 

  Для произвольной числовой полугруппы $S$ и $d\in\mathbb Z_{>0}$ в \cite{1} доказана формула тригонометрического типа для рода числовой полугруппы
$$
S/d:=\{s\in\mathbb Z_{\ge 0}:ds\in S\}.
$$
В п.6 мы докажем формулу рационального типа для рода числовой полугруппы $S/d$.

\bigskip

\section{Числовые полугруппы и их множества Апери}

\bigskip

  Множество $S\subset \mathbb Z_{\ge 0}$ называется числовой полугруппой, если оно содержит $0$, замкнуто относительно сложения и имеет конечное дополнение в $\mathbb Z_{\ge 0}$. Для ненулевого $s\in S$ и $I_s:=\{0,1,\dots,s-1\}$ множество Апери $Ap_s(S)$ определяется следующим образом:
$$
Ap_s(S):=\{u\in S:u-s\notin S\}
$$
$$
=\{a_0,a_1,\dots,a_{s-1}:a_k=\min((k+s\mathbb Z_{\ge 0})\cap S),k\in I_s\}.
$$
Другими словами, $Ap_s(S)$ состоит из минимальных элементов в пересечениях полугруппы $S$ с классами вычетов $(k+s\mathbb Z_{\ge 0})$ -- для каждого $k\in I_s$
$$
\{k,k+s,\dots,k+s([a_k/s]-1)=a_k-s\}\subset C_{S}:=\mathbb Z_{\ge 0}\setminus S,      
$$
если $[a_k/s]>0$ и
$$
(k+s\mathbb Z_{\ge 0})\cap C_{S}=\emptyset,      
$$
если $[a_k/s]=0$ (например, $a_0=0$). Пусть $C_{S}(q):=\sum_{k\in C_{S}}q^k$.

\begin{prop}\label{prop1} Для ненулевого $s\in S$ и каждого $k\in I_s$
$$
q^{s\left[\frac{a_k}{s}\right]}=1+\frac{q^s-1}{sq^k}\sum_{j=0}^{s-1}\epsilon_s^{-jk}C_{S}(\epsilon_s^jq),           \eqno (2)
$$
$$
\left[\frac{a_k}{s}\right]=\frac{1}{s}\sum_{j=0}^{s-1}\epsilon_s^{-jk}C_{S}(\epsilon_s^j),        \eqno (3)
$$
в частности, для каждого $k\in I_b$
$$
q^{b\left[\frac{ak}{b}\right]}=1+\frac{q^b-1}{bq^{\pi(k)}}\sum_{j=0}^{b-1}\epsilon_b^{-jak}C_{S_{a,b}}(\epsilon_b^jq),   \eqno (4)
$$
$$
\left[\frac{ak}{b}\right]=\frac{1}{b}\sum_{j=0}^{b-1}\epsilon_b^{-jak}C_{S_{a,b}}(\epsilon_b^j).     \eqno (5)
$$
\end{prop}

Доказательство. 
$$
C_{S}(q)=\sum_{k=1}^{s-1}\sum_{i\in(k+s\mathbb Z_{\ge 0})\cap C_S}q^i=\sum_{k=1}^{s-1}(q^k+q^{s+k}+\dots+q^{s([a_k/s]-1)+k})
$$
$$
=\sum_{k=1}^{s-1}q^k\frac{q^{s[a_k/s]}-1}{q^s-1}=\sum_{k=1}^{s-1}\frac{q^{a_k}-q^k}{q^s-1} \eqno (6)
$$
Это равенство можно обратить, а именно, хорошо известен способ выделить слагаемые степенного ряда, степени которых принадлежат фиксированному классу вычетов по модулю $s$ (multisection series method),
$$
q^k\frac{q^{s[a_k/s]}-1}{q^s-1}=\frac{1}{s}\sum_{j=0}^{s-1}\epsilon_s^{-jk}C_{S}(\epsilon_s^jq), \qquad  k\in I_s.
$$
Отсюда получаем формулу (2) и, переходя к пределу при $q\to1$, формулу (3). Формулы (4) и (5) следуют из известного равенства
$$
Ap_b(S_{a,b})=\{0a,1a,\dots,(b-1)a\}.        
$$
По определению числа $k$ и $a_k\in Ap_s(S)$ имеют остаток $k$ при делении на $s$. Для $Ap_b(S_{a,b})$ такими парами являются $\pi(k)$ и $ak$. Поэтому
$$
q^{b\left[\frac{ak}{b}\right]}=1+\frac{q^b-1}{bq^{\pi(k)}}\sum_{j=0}^{b-1}\epsilon_b^{-j\pi(k)}C_{S_{a,b}}(\epsilon_b^jq),   
$$
$$
\left[\frac{ak}{b}\right]=\frac{1}{b}\sum_{j=0}^{b-1}\epsilon_b^{-j\pi(k)}C_{S_{a,b}}(\epsilon_b^j)     
$$
и $\epsilon_b^{-j\pi(k)}=\epsilon_b^{-jak}$.

\bigskip

\section{Суммы Вороного}

\bigskip

  Многочлены Апостола-Бернулли $B_k(q,\lambda)$ определяются производящей функцией
$$
\frac{te^{tq}}{\lambda e^t-1}=\sum_{k=0}^\infty B_k(q,\lambda)\frac{t^k}{k!}.
$$
Они связаны с многочленами Мириманова 
$$
M_{b-1}(\lambda,m):=\sum_{k=0}^{b-1}k^m\lambda^k
$$
формулой \cite{3}, \cite{17}
$$
M_{b-1}(\lambda,m)=\frac{\lambda^bB_{m+1}(b,\lambda)-B_{m+1}(0,\lambda)}{m+1}.
$$

\begin{prop}\label{prop2} Для $m>0$
$$
V_{m,n}(a,b)=\frac{1}{b^n}\sum_{\substack{0\le i_0,\dots,i_{b-1}\le n\\i_0+\dots+i_{b-1}=n}}\binom{n}{i_0\dots i_{b-1}}\prod_{j=0}^{b-1}C_{S_{a,b}}^{i_j}(\epsilon_b^j)
$$
$$
\times M_{b-1}(\epsilon_b^{-a\sum_{j=0}^{b-1}ji_j},m)
$$
$$
=\frac{1}{b^n}\sum_{\substack{0\le i_0,\dots,i_{b-1}\le n\\i_0+\dots+i_{b-1}=n}}\binom{n}{i_0\dots i_{b-1}}\prod_{j=0}^{b-1}C_{S_{a,b}}^{i_j}(\epsilon_b^j)
$$
$$
\times\frac{B_{m+1}(b,\epsilon_b^{-a\sum_{j=0}^{b-1}ji_j})-B_{m+1}(0,\epsilon_b^{-a\sum_{j=0}^{b-1}ji_j})}{m+1},
$$
в частности,
$$
V_{m,1}(a,b)=\frac{1}{b}\sum_{j=0}^{b-1}C_{S_{a,b}}(\epsilon_b^j)M_{b-1}(\epsilon_b^{-aj},m)
$$
$$
=\frac{1}{b}\sum_{j=0}^{b-1}C_{S_{a,b}}(\epsilon_b^j)\frac{B_{m+1}(b,\epsilon_b^{-aj})-B_{m+1}(0,\epsilon_b^{-aj})}{m+1}.
$$
\end{prop}

Доказательство. 
$$
V_{m,n}(a,b)=\sum_{k=1}^{b-1}k^m\left[\frac{ak}{b}\right]^n=\frac{1}{b^n}\sum_{k=1}^{b-1}k^m\left(\sum_{j=0}^{b-1}\epsilon_b^{-jak}C_{S_{a,b}}(\epsilon_b^j)\right)^n
$$
$$
=\frac{1}{b^n}\sum_{k=1}^{b-1}k^m\sum_{\substack{0\le i_0,\dots,i_{b-1}\le n\\i_0+\dots+i_{b-1}=n}}\binom{n}{i_0\dots i_{b-1}}\epsilon_b^{-ak\sum_{j=0}^{b-1}ji_j}\prod_{j=0}^{b-1}C_{S_{a,b}}^{i_j}(\epsilon_b^j)
$$
$$
=\frac{1}{b^n}\sum_{\substack{0\le i_0,\dots,i_{b-1}\le n\\i_0+\dots+i_{b-1}=n}}\binom{n}{i_0\dots i_{b-1}}\prod_{j=0}^{b-1}C_{S_{a,b}}^{i_j}(\epsilon_b^j)\sum_{k=1}^{b-1}k^m\epsilon_b^{-ak\sum_{j=0}^{b-1}ji_j}
$$

\bigskip

\section{Многочлены Дедекинда-Карлитца}

\bigskip

  Перестановка Золотарёва на множестве $I_b$ определяется равенством
$$
\pi(k)=ak\mod b, \qquad k\in I_b.
$$
Другими словами, деление с остатком числа $ak$ на число $b$ запишется в следующем виде
$$
ak=b\left[\frac{ak}{b}\right]+\pi(k).
$$
Эта перестановка использовалось Е. И. Золотарёвым в 1872 г. при доказательстве квадратичного закона взаимности.

Введём обозначение
$$
d_j(q,t;a,b):=\sum_{k=1}^{b-1}\frac{\epsilon_b^{-jak}t^{k-1}}{q^{\pi(k)}}.
$$

\begin{prop}\label{prop3} 
$$
c(q^b,t;a,b):=\frac{t^{b-1}-1}{t-1}+\frac{q^b-1}{b}\sum_{j=0}^{b-1}d_j(q,t;a,b)C_{S_{a,b}}(\epsilon_b^jq).
$$
\end{prop}

Доказательство. Используя формулу (4), получим
$$
c(q^b,t;a,b):=\sum_{k=1}^{b-1}q^{b[ak/b]}t^{k-1}
$$
$$
=\sum_{k=1}^{b-1}\left(1+\frac{q^b-1}{bq^{\pi(k)}}\sum_{j=0}^{b-1}\epsilon_b^{-jak}C_{S_{a,b}}(\epsilon_b^jq)\right)t^{k-1}
$$
$$
=\frac{t^{b-1}-1}{t-1}+\frac{q^b-1}{b}\sum_{j=0}^{b-1}C_{S_{a,b}}(\epsilon_b^jq)\sum_{k=1}^{b-1}\frac{\epsilon_b^{-jak}t^{k-1}}{q^{\pi(k)}}.
$$

  Определим многочлены
$$
R_{m,n}(q,t;a,b):=\sum_{k=1}^{b-1}\left(\frac{q^{[ak/b]}-1}{q-1}\right)^n\left(\frac{t^k-1}{t-1}\right)^m,
$$
$$
T_{m,n}(q,t;a,b):=\sum_{k=1}^{b-1}\left(\frac{q^{ak}-q^{\pi(k)}}{q^b-1}\right)^n\left(\frac{t^k-1}{t-1}\right)^m.
$$
Тогда
$$
R_{m,n}(1,1;a,b)=T_{m,n}(1,1;a,b)=V_{m,n}(a,b).
$$

\begin{prop}\label{prop4} 
$$
R_{1,1}(q,t;a,b)=\frac{tc(q,t;a,b)-(b-1)(q^{a-1}-1)-\frac{t^b-t}{t-1}}{(q-1)(t-1)}
$$
$$
+\frac{1}{t-1}\sum_{k=1}^{a-1}[bk/a]q^{k-1},
$$
$$
T_{1,1}(q,t;a,b)=q^{\pi(k)}\frac{tc(q^b,t;a,b)-(b-1)(q^{b(a-1)}-1)-\frac{t^b-t}{t-1}}{(q^b-1)(t-1)}
$$
$$
+\frac{q^{\pi(k)}}{t-1}\sum_{k=1}^{a-1}[bk/a]q^{b(k-1)}.
$$
\end{prop}

Доказательство.
$$
R_{1,1}(q,t;a,b)=\frac{\sum_{k=1}^{b-1}(q^{[ak/b]}t^k-q^{[ak/b]}-t^k+1)}{(q-1)(t-1)},
$$
$$
T_{1,1}(q,t;a,b)=q^{\pi(k)}\frac{\sum_{k=1}^{b-1}(q^{b[ak/b]}t^k-q^{b[ak/b]}-t^k+1)}{(q^b-1)(t-1)}
$$
и далее используем формулу из \cite{6}, Corollary 5.10,
$$
\sum_{k=1}^{b-1}q^{[ak/b]}=q^{a-1}(b-1)-(q-1)\sum_{k=1}^{a-1}[bk/a]q^{k-1}.
$$

  Многочлен
$$
\sum_{k=0}^{b-1}\left(\frac{ak}{b}-\left[\frac{ak}{b}\right]-\frac{1}{2}\right)q^k=\frac{a}{b}\frac{bq^b(q-1)-q(q^b-1)}{(q-1)^2}-\frac{q^b-1}{2(q-1)}
$$
$$
-\sum_{k=0}^{b-1}\left[\frac{ak}{b}\right]q^k
$$
появляется в \cite{8}, где для него доказана теорема взаимности.

\begin{prop}\label{prop4} 
$$
\sum_{k=0}^{b-1}\left[\frac{ak}{b}\right]q^k=\frac{q^b-1}{b}\sum_{j=0}^{b-1}\frac{C_{S_{a,b}}(\epsilon_b^j)}{\epsilon_b^{-ja}q-1}.
$$
\end{prop}

Доказательство. Используя формулу (5), получим
$$
\sum_{k=0}^{b-1}\left[\frac{ak}{b}\right]q^k=\frac{1}{b}\sum_{k=0}^{b-1}q^k\sum_{j=0}^{b-1}\epsilon_b^{-jak}C_{S_{a,b}}(\epsilon_b^j)
$$
$$
=\frac{1}{b}\sum_{j=0}^{b-1}C_{S_{a,b}}(\epsilon_b^j)\sum_{k=0}^{b-1}(\epsilon_b^{-ja}q)^k.
$$

\bigskip

\section{Классические суммы Дедекинда}

\bigskip

  Классические суммы Дедекинда можно определить одним из следующих равенств
$$
s(a,b):=\sum_{k=1}^{b-1}\left(\left(\frac{k}{b}\right)\right)\left(\left(\frac{ak}{b}\right)\right)
$$
$$
=\sum_{k=1}^{b-1}\frac{k}{b}\left(\left(\frac{ak}{b}\right)\right)
$$
$$
=-\frac{1}{b}\sum_{k=1}^{b-1}\frac{1}{(\epsilon_b^{ak}-1)(\epsilon_b^k-1)}+\frac{b-1}{4b}
$$
$$
=\frac{1}{4b}\sum_{k=1}^{b-1}\frac{1+\epsilon_b^k}{1-\epsilon_b^k}\frac{1+\epsilon_b^{-ak}}{1-\epsilon_b^{-ak}}
$$
$$
=\frac{1}{4b}\sum_{k=1}^{b-1}\cot\frac{\pi k}{b}\cot\frac{\pi ak}{b}
$$
$$
=-\frac{1}{b}V_{1,1}(a,b)+\frac{b-1}{4}\left(\frac{4a}{3}-\frac{2a}{3b}-1\right),
$$
где $((x))=x-[x]-1/2$, если $x\notin \mathbb Z$, и $((x))=0$, если $x\in \mathbb Z$ \cite{13}.

  По-видимому, из этих равенств легко вывести формулу (7), учитывая явный вид чисел $C_{S_{a,b}}(\epsilon_b^j)$: для $0<k<b$
$$
C_{S_{a,b}}(\epsilon_b^k)=\left.\frac{A_{S_{a,b}}(q)-1}{q-1}\right|_{q\to\epsilon_b^k}
$$
$$
=\left.\left(\frac{1}{1-q}-H_{S_{a,b}}(q)\right)\right|_{q\to\epsilon_b^k}=\frac{a}{\epsilon_b^{ka}-1}-\frac{1}{\epsilon_b^k-1},
$$
число
$$
C_{S_{a,b}}(1)=|C_{S_{a,b}}|=\left.\frac{A_{S_{a,b}}(q)-1}{q-1}\right|_{q\to1}
$$
$$
=\left.\left(\frac{1}{1-q}-H_{S_{a,b}}(q)\right)\right|_{q\to1}=\frac{(a-1)(b-1)}{2}
$$
известно как род числовой полугруппы $S_{a,b}$ (или $(a,b)$-торического узла). Но в данном случае нас интересует приходящее из числовой полугруппы $S_{a,b}$ доказательство.

\begin{prop}\label{prop5}
$$
V_{1,1}(a,b)=\sum_{k=1}^{b-1}k\left[\frac{ak}{b}\right]=\sum_{j=1}^{b-1}\frac{C_{S_{a,b}}(\epsilon_b^j)}{\epsilon_b^{-ja}-1}+\frac{(a-1)(b-1)^2}{4}.         \eqno (7)
$$
\end{prop}

Доказательство. Используя формулу (5), получим
$$
V_{1,1}(a,b)=\frac{1}{b}\sum_{j=0}^{b-1}C_{S_{a,b}}(\epsilon_b^j)\sum_{k=1}^{b-1}k\epsilon_b^{-jak}.
$$
Упростим коэффициенты при $C_{S_{a,b}}(\epsilon_b^j)$. Для $0<j<b$
$$
\frac{1}{b}\sum_{k=1}^{b-1}k\epsilon_b^{-jak}=\frac{1}{b}\left.\left(q\frac{d}{dq}\sum_{k=0}^{b-1}q^k\right)\right|_{q=\epsilon_b^{-ja}}
=\frac{1}{b}\left.\left(q\frac{d}{dq}\frac{q^b-1}{q-1}\right)\right|_{q=\epsilon_b^{-ja}}
$$
$$
=\frac{1}{b}\left.\frac{bq^b(q-1)-q(q^b-1)}{(q-1)^2}\right|_{q=\epsilon_b^{-ja}}
=\frac{1}{\epsilon_b^{-ja}-1}.
$$

\bigskip

\section{Род числовой полугруппы $S/d$}

\bigskip

  Для $d\in\mathbb Z_{>0}$ формула для рода числовой полугруппы $S/d$ из \cite{1} имеет вид
$$
g(S/d)=\frac{1}{d}\left(g(S)+\sum_{k=1}^{d-1}C_S(\epsilon_d^k)\right)
$$
$$
=\frac{1}{d}\left(g(S)+\sum_{k=1}^{d-1}\frac{A_S(\epsilon_d^k)-1}{\epsilon_d^k-1}\right)=\frac{1}{d}\left(g(S)+\frac{d-1}{2}-\sum_{k=1}^{d-1}H_S(\epsilon_d^k)\right).
$$
Для ненулевого $s\in S/d$ пусть
$$
Ap_{ds}(S)=\{a_0,a_1,\dots,a_{ds-1}:a_l=l\mod ds,l\in I_{ds}\}.
$$

\begin{prop}\label{prop6} Для ненулевого $s\in S/d$ и $d\in\mathbb Z_{>0}$
$$
g(S/d)=\sum_{i=1}^{s-1}\left[\frac{a_{di}}{ds}\right].
$$
\end{prop}

Доказательство. Переходя в равенстве (6) для $Ap_{ds}(S)$ к пределу при $q\to\epsilon_{ds}^k$, получим для $k\in I_{ds}$
$$
C_{S}(\epsilon_{ds}^k)=\sum_{l=1}^{ds-1}\epsilon_{ds}^{kl}\left[\frac{a_l}{ds}\right],
$$
в частности,
$$
g(S)=C_{S}(1)=\sum_{l=1}^{ds-1}\left[\frac{a_l}{ds}\right].
$$
Поэтому
$$
g(S/d)=\frac{1}{d}\left(g(S)+\sum_{k=1}^{d-1}C_S(\epsilon_d^k)\right)=\frac{1}{d}\sum_{k=0}^{d-1}C_S(\epsilon_{ds}^{ks})
$$
$$
=\frac{1}{d}\sum_{k=0}^{d-1}\sum_{l=1}^{ds-1}\epsilon_{ds}^{kls}\left[\frac{a_l}{ds}\right]=\frac{1}{d}\sum_{l=1}^{ds-1}\left[\frac{a_l}{ds}\right]\sum_{k=0}^{d-1}\epsilon_{d}^{kl}
=\sum_{i=1}^{s-1}\left[\frac{a_{di}}{ds}\right].
$$


\begin{thebibliography}{24}

\bibitem{1}

A. Adeniran, S. Butler, C. Defant, Y. Gao, P.E. Harris, C. Hettle, Q. Liang, H. Nam, A. Volk, On the genus of a quotient of a numerical semigroup, Semigroup Forum 98 (2019), 690-700.

\bibitem{2}

T. Agoh, Voronoi type congruence and its applications, Eur. J. Pure Appl. Math. 1 (2008),
3-21.

\bibitem{3}

T. M. Apostol, Addendum to “On the Lerch Zeta-function,” Pacific J. Math. 2 (1952), 10.

\bibitem{4}

L. Bardomero, M. Beck, Frobenius coin-exchange generating functions, Amer. Math. Monthly 127 (2020), no. 4, 308-315.

\bibitem{5}

M. Beck, C. Haase, A. R. Matthews, Dedekind-Carlitz polynomials as lattice-point
enumerators in rational polyhedra, Math. Ann. 341 (2008), 945-961.

\bibitem{6}

B. C. Berndt, U. Dieter, Sums involving the greatest integer function and Riemann-Stieltjes integration, J. Reine Angew. Math. 337 (1982), 208-220.

\bibitem{7}

L. Carlitz, The number of terms in the cyclotomic polynomial $F_{pq}(x)$, Amer. Math. Monthly 73 (1966), 979-981.

\bibitem{8}

L. Carlitz, Some polynomials associated with Dedekind sums, Acta Math. Acad. Sci. Hungar. 26(3-4) (1975), 311-319.

\bibitem{9}

L. Carlitz, Some sums containing the greatest integer function, Fibonacci Quar. 15(1) (1977), 78-84.

\bibitem{10}

L. Carlitz, Some restricted multiple sums, Fibonacci Quar. 18(1) (1980), 58-65.

\bibitem{11}

L.J. Mordell, The reciprocity formula for Dedekind sums, Amer. J. Math. 73 (1951), 593-598.

\bibitem{12}

P. Moree, Numerical semigroups, cyclotomic polynomials, and Bernoulli numbers, Amer. Math. Monthly 121(10) (2014), 890-902.

\bibitem{13}

H. Rademacher, E. Grosswald, Dedekind Sums, The Carus Mathematical Monographs, No. 16 (Mathematical Association of America, 1972).

\bibitem{14}

S. Sertoz, A. E. Ozluk, On the number of representations of an integer by a linear form, Istanb. Univ. Fen Fak. Mat. Derg. 50 (1991), 67-77.

\bibitem{15}

L. A. Szekely, N. C. Wormald, Generating functions for the Frobenius Problem with 2 and 3 generators, Mathematical Chronicle 15 (1986), 49-57.

\bibitem{16}

J. V. Uspensky, M. A. Heaslet, Elementary Number Theory, McGraw-Hill, New York, 1939.

\bibitem{17}

W.-P. Wang, W.-W. Wang, Some results on power sums and Apostol-type polynomials, Integral Transforms Spec. Funct. 21 (2010), 307-318.

\end{thebibliography}
\end {document}